\documentclass{amsart}
\usepackage{amssymb}
\usepackage{amsmath}
\usepackage{amsfonts}

\newtheorem{theorem}{Theorem}
\theoremstyle{plain}

\newtheorem{corollary}{Corollary}

\newtheorem{definition}{Definition}

\newtheorem{proposition}{Proposition}

\numberwithin{equation}{section}
\input{tcilatex}

\begin{document}
\title[$q$\textbf{-Hardy-littlewood-type maximal operator with weight}]{%
\textbf{On the} $q$\textbf{-Hardy-littlewood-type maximal operator with
weight related to fermionic }$p$\textbf{-adic }$q$\textbf{-integral on }$%
\mathbb{Z}
_{p}$}
\author{\textbf{Serkan Araci}}
\address{\textbf{University of Gaziantep, Faculty of Science and Arts,
Department of Mathematics, 27310 Gaziantep, TURKEY}}
\email{\textbf{mtsrkn@hotmail.com}}
\author{\textbf{Mehmet Acikgoz}}
\address{\textbf{University of Gaziantep, Faculty of Science and Arts,
Department of Mathematics, 27310 Gaziantep, TURKEY}}
\email{\textbf{acikgoz@gantep.edu.tr}}

\begin{abstract}
The fundamental aim of this paper is to define weighted $q$%
-Hardy-littlewood-type maximal operator by means of fermionic $p$-adic $q$%
-invariant distribution on $%
\mathbb{Z}
_{p}$. Also, we derive some interesting properties concerning this type
maximal operator.
\end{abstract}

\maketitle

\section{\textbf{Introduction and Notations}}

\bigskip $p$-adic numbers also play a vital and important role in
mathematics. $p$-adic numbers were invented by the German mathematician Kurt
Hensel \cite{Hensel}, around the end of the nineteenth century.\ In spite of
their being already one hundred years old, these numbers are still today
enveloped in an aura of mystery within the scientific community.

The fermionic $p$-adic $q$-integral are originally constructed by Kim \cite%
{Kim 4}. Kim also introduced Lebesgue-Radon-Nikodym Theorem with respect to
fermionic $p$-adic $q$-integral on $%
\mathbb{Z}
_{p}$. The fermionic $p$--adic $q$-integral on $%
\mathbb{Z}
_{p}$ is used in Mathematical Physics for example the functional equation of
the $q$-Zeta function, the $q$-Stirling numbers and $q$-Mahler theory of
integration with respect to the ring $%
\mathbb{Z}
_{p}$ together with Iwasawa's $p$-adic $q$-$L$ function.

In \cite{Jang}, Jang also defined $q$-extension of Hardy-Littlewood-type
maximal operator by means of $q$-Volkenborn integral on $%
\mathbb{Z}
_{p}$. Next, in previous paper \cite{Araci}, Araci and Acikgoz added a
weight into Jang's $q$-Hardy-Littlewood-type maximal operator and derived
some interesting properties by means of Kim's $p$-adic $q$-integral on $%
\mathbb{Z}
_{p}$. Now also, we shall consider weighted $q$-Hardy-Littlewood-type
maximal operator on the fermionic $p$-adic $q$-integral on $%
\mathbb{Z}
_{p}$. Moreover, we shall analyse $q$-Hardy-Littlewood-type maximal operator
via the fermionic $p$-adic $q$-integral on $%
\mathbb{Z}
_{p}$.

Assume that $p$ be an odd prime number. Let $\mathcal{%
\mathbb{Q}
}_{p}$ be the field of $p$-adic rational numbers and let $\mathcal{%
\mathbb{C}
}_{p}$ be the completion of algebraic closure of $\mathcal{%
\mathbb{Q}
}_{p}$.

Thus, 
\begin{equation*}
\mathcal{%
\mathbb{Q}
}_{p}=\left\{ x=\sum_{n=-k}^{\infty }a_{n}p^{n}:0\leq a_{n}<p\right\} .
\end{equation*}

Then $%
\mathbb{Z}
_{p}$ is an integral domain, which is defined by 
\begin{equation*}
\mathcal{%
\mathbb{Z}
}_{p}=\left\{ x=\sum_{n=0}^{\infty }a_{n}p^{n}:0\leq a_{n}\leq p-1\right\} ,
\end{equation*}

or 
\begin{equation*}
\mathcal{%
\mathbb{Z}
}_{p}=\left\{ x\in 
\mathbb{Q}
_{p}:\left\vert x\right\vert _{p}\leq 1\right\} .
\end{equation*}

In this paper, we assume that $q\in 
\mathbb{C}
_{p}$ with $\left\vert 1-q\right\vert _{p}<1$ as an indeterminate.

The $p$-adic absolute value $\left\vert .\right\vert _{p}$, is normally
defined by 
\begin{equation*}
\left\vert x\right\vert _{p}=\frac{1}{p^{r}}\text{,}
\end{equation*}

where $x=p^{r}\frac{s}{t}$ with $\left( p,s\right) =\left( p,t\right)
=\left( s,t\right) =1$ and $r\in \mathcal{%
\mathbb{Q}
}$.

A $p$-adic Banach space $B$ is a $%
\mathbb{Q}
_{p}$-vector space with a lattice $B^{0}$ ($\mathcal{%
\mathbb{Z}
}_{p}$-module) separated and complete for $p$-adic topology, ie., 
\begin{equation*}
B^{0}\simeq \lim_{\overleftarrow{n\in 
\mathbb{N}
}}B^{0}/p^{n}B^{0}\text{.}
\end{equation*}

For all $x\in B$, there exists $n\in \mathcal{%
\mathbb{Z}
}$, such that $x\in p^{n}B^{0}$. Define 
\begin{equation*}
v_{B}\left( x\right) =\sup_{n\in 
\mathbb{N}
\cup \left\{ +\infty \right\} }\left\{ n:x\in p^{n}B^{0}\right\} \text{.}
\end{equation*}

It satisfies the following properties:%
\begin{eqnarray*}
v_{B}\left( x+y\right) &\geq &\min \left( v_{B}\left( x\right) ,v_{B}\left(
y\right) \right) \text{,} \\
v_{B}\left( \beta x\right) &=&v_{p}\left( \beta \right) +v_{B}\left(
x\right) \text{, if }\beta \in 
\mathbb{Q}
_{p}\text{.}
\end{eqnarray*}

Then, $\left\Vert x\right\Vert _{B}=p^{-v_{B}\left( x\right) }$ defines a
norm on $B,$ such that $B$ is complete for $\left\Vert .\right\Vert _{B}$
and $B^{0}$ is the unit ball.

A measure on $\mathcal{%
\mathbb{Z}
}_{p}$ with values in a $p$-adic Banach space $B$ is a continuous linear map%
\begin{equation*}
f\mapsto \int f\left( x\right) \mu =\int_{%
\mathbb{Z}
_{p}}f\left( x\right) \mu \left( x\right)
\end{equation*}

from $C^{0}\left( \mathcal{%
\mathbb{Z}
}_{p},\mathcal{%
\mathbb{C}
}_{p}\right) $, (continuous function on $\mathcal{%
\mathbb{Z}
}_{p}$) to $B$. We know that the set of locally constant functions from $%
\mathcal{%
\mathbb{Z}
}_{p}$ to $\mathcal{%
\mathbb{Q}
}_{p}$ is dense in $C^{0}\left( \mathcal{%
\mathbb{Z}
}_{p},\mathcal{%
\mathbb{C}
}_{p}\right) $ so.

Explicitly, for all $f\in C^{0}\left( \mathcal{%
\mathbb{Z}
}_{p},\mathcal{%
\mathbb{C}
}_{p}\right) $, the locally constant functions 
\begin{equation*}
f_{n}=\sum_{i=0}^{p^{n}-1}f\left( i\right) 1_{i+p^{n}%
\mathbb{Z}
_{p}}\rightarrow \text{ }f\text{ in }C^{0}\text{.}
\end{equation*}

Now if ~$\mu \in \mathcal{D}_{0}\left( \mathcal{%
\mathbb{Z}
}_{p},\mathcal{%
\mathbb{Q}
}_{p}\right) $, set $\mu \left( i+p^{n}\mathcal{%
\mathbb{Z}
}_{p}\right) =\int_{%
\mathbb{Z}
_{p}}1_{i+p^{n}\mathcal{%
\mathbb{Z}
}_{p}}\mu $. Then $\int_{\mathcal{%
\mathbb{Z}
}_{p}}f\mu $ is given by the following \textquotedblleft Riemann
sums\textquotedblright 
\begin{equation*}
\int_{%
\mathbb{Z}
_{p}}f\mu =\lim_{n\rightarrow \infty }\sum_{i=0}^{p^{n}-1}f\left( i\right)
\mu \left( i+p^{n}\mathcal{%
\mathbb{Z}
}_{p}\right) \text{.}
\end{equation*}

T. Kim defined $\mu _{-q}$ as follows:%
\begin{equation*}
\mu _{-q}\left( \xi +dp^{n}\mathcal{%
\mathbb{Z}
}_{p}\right) =\frac{\left( -q\right) ^{\xi }}{\left[ dp^{n}\right] _{-q}}
\end{equation*}

and this can be extended to a distribution on $\mathcal{%
\mathbb{Z}
}_{p}$. This distribution yields an integral in the case $d=1$.

So, $q$-Volkenborn integral was defined by T. Kim as follows:%
\begin{equation}
I_{-q}\left( f\right) =\int_{\mathcal{%
\mathbb{Z}
}_{p}}f\left( \xi \right) d\mu _{q}\left( \xi \right) =\lim_{n\rightarrow
\infty }\frac{1}{\left[ p^{n}\right] _{-q}}\sum_{\xi =0}^{p^{n}-1}\left(
-1\right) ^{\xi }f\left( \xi \right) q^{\xi }\text{ }  \label{equation 6}
\end{equation}

Where $\left[ x\right] _{q}$ is a $q$-extension of $x$ which is defined by%
\begin{equation*}
\left[ x\right] _{q}=\frac{1-q^{x}}{1-q}\text{,}
\end{equation*}

note that $\lim_{q\rightarrow 1}\left[ x\right] _{q}=x$ cf. \cite{Kim 2}, 
\cite{Kim 3}, \cite{Kim 4}, \cite{Kim 5}, \cite{Jang}.

Let $d$ be a fixed positive integer with $\left( p,d\right) =1$. We now set%
\begin{eqnarray*}
X &=&X_{d}=\lim_{\overleftarrow{n}}\mathcal{%
\mathbb{Z}
}/dp^{n}\mathcal{%
\mathbb{Z}
}, \\
X_{1} &=&%
\mathbb{Z}
_{p}, \\
X^{\ast } &=&\underset{\underset{\left( a,p\right) =1}{0<a<dp}}{\cup }a+dp%
\mathcal{%
\mathbb{Z}
}_{p}, \\
a+dp^{n}\mathcal{%
\mathbb{Z}
}_{p} &=&\left\{ x\in X\mid x\equiv a\left( \func{mod}p^{n}\right) \right\} ,
\end{eqnarray*}

where $a\in \mathcal{%
\mathbb{Z}
}$ satisfies the condition $0\leq a<dp^{n}$. For $f\in UD\left( \mathcal{%
\mathbb{Z}
}_{p},\mathcal{%
\mathbb{C}
}_{p}\right) $,%
\begin{equation*}
\int_{%
\mathbb{Z}
_{p}}f\left( x\right) d\mu _{-q}\left( x\right) =\int_{X}f\left( x\right)
d\mu _{-q}\left( x\right) ,
\end{equation*}

(for details, see \cite{Kim 8}).

By the meaning of $q$-Volkenborn integral, we consider below strongly $p$%
-adic $q$-invariant distribution $\mu _{-q}$ on $%
\mathbb{Z}
_{p}$ in the form 
\begin{equation*}
\left\vert \left[ p^{n}\right] _{-q}\mu _{-q}\left( a+p^{n}\mathcal{%
\mathbb{Z}
}_{p}\right) -\left[ p^{n+1}\right] _{-q}\mu _{-q}\left( a+p^{n+1}\mathcal{%
\mathbb{Z}
}_{p}\right) \right\vert <\delta _{n},
\end{equation*}

where $\delta _{n}\rightarrow 0$ as $n\rightarrow \infty $ and $\delta _{n}$
is independent of $a$. Let $f\in UD\left( \mathcal{%
\mathbb{Z}
}_{p},\mathcal{%
\mathbb{C}
}_{p}\right) $, for any $a\in \mathcal{%
\mathbb{Z}
}_{p}$, we assume that the weight function $\omega \left( x\right) $ is
defined by $\omega \left( x\right) =\omega ^{x}$ where $\omega \in 
\mathbb{C}
_{p}$ with $\left\vert 1-\omega \right\vert _{p}<1$. We define the weighted
measure on $\mathcal{%
\mathbb{Z}
}_{p}$ as follows:%
\begin{equation}
\mu _{f,-q}^{\left( \omega \right) }\left( a+p^{n}\mathcal{%
\mathbb{Z}
}_{p}\right) =\int_{a+p^{n}\mathcal{%
\mathbb{Z}
}_{p}}\omega ^{\xi }f\left( \xi \right) d\mu _{-q}\left( \xi \right)
\label{equation 2}
\end{equation}

where the integral is the fermionic $p$-adic $q$-integral. By (\ref{equation
2}), we easily note that $\mu _{f,-q}^{\left( \omega \right) }$ is a
strongly weighted measure on $%
\mathbb{Z}
_{p}$. Namely,%
\begin{eqnarray*}
&&\left\vert \left[ p^{n}\right] _{-q}\mu _{f,-q}^{\left( \omega \right)
}\left( a+p^{n}\mathcal{%
\mathbb{Z}
}_{p}\right) -\left[ p^{n+1}\right] _{-q}\mu _{f,-q}^{\left( \omega \right)
}\left( a+p^{n+1}\mathcal{%
\mathbb{Z}
}_{p}\right) \right\vert _{p} \\
&=&\left\vert \sum_{x=0}^{p^{n}-1}\left( -1\right) ^{x}\omega ^{x}f\left(
x\right) q^{x}-\sum_{x=0}^{p^{n}}\left( -1\right) ^{x}\omega ^{x}f\left(
x\right) q^{x}\right\vert _{p} \\
&\leq &\left\vert \frac{f\left( p^{n}\right) \left( -1\right) ^{p^{n}}\omega
^{p^{n}}q^{p^{n}}}{p^{n}}\right\vert _{p}\left\vert p^{n}\right\vert _{p} \\
&\leq &Cp^{-n}
\end{eqnarray*}

Thus, we get the following proposition.

\begin{proposition}
For $f,g\in UD\left( \mathcal{%
\mathbb{Z}
}_{p},\mathcal{%
\mathbb{C}
}_{p}\right) $, then, we have 
\begin{equation*}
\mu _{\alpha f+\beta g,-q}^{\left( \omega \right) }\left( a+p^{n}\mathcal{%
\mathbb{Z}
}_{p}\right) =\alpha \mu _{f,-q}^{\left( \omega \right) }\left( a+p^{n}%
\mathcal{%
\mathbb{Z}
}_{p}\right) +\beta \mu _{g,-q}^{\left( \omega \right) }\left( a+p^{n}%
\mathcal{%
\mathbb{Z}
}_{p}\right) \text{.}
\end{equation*}%
where $\alpha ,\beta $ are positive constants. Also, we have%
\begin{equation*}
\left\vert \left[ p^{n}\right] _{-q}\mu _{f,-q}^{\left( \omega \right)
}\left( a+p^{n}\mathcal{%
\mathbb{Z}
}_{p}\right) -\left[ p^{n+1}\right] _{-q}\mu _{f,-q}^{\left( \omega \right)
}\left( a+p^{n+1}\mathcal{%
\mathbb{Z}
}_{p}\right) \right\vert \leq Cp^{-n}
\end{equation*}%
where $C$ is positive constant.
\end{proposition}

Let $\mathcal{P}_{q}\left( x\right) \in 
\mathbb{C}
_{p}\left[ \left[ x\right] _{q}\right] $ be an arbitrary $q$-polynomial. Now
also, we indicate that $\mu _{\mathcal{P},-q}^{\left( \omega \right) }$ is a
strongly weighted fermionic $p$-adic $q$-invariant measure on $%
\mathbb{Z}
_{p}$. Without a loss of generality, it is sufficient to evidence the
statement for $\mathcal{P}\left( x\right) =\left[ x\right] _{q}^{k}$. 
\begin{equation}
\mu _{\mathcal{P},-q}^{\left( \omega \right) }\left( a+p^{n}\mathcal{%
\mathbb{Z}
}_{p}\right) =\lim_{m\rightarrow \infty }\frac{1}{\left[ p^{m}\right] _{-q}}%
\sum_{i=0}^{p^{m-n}-1}w^{a+ip^{n}}\left[ a+ip^{n}\right] _{q}^{k}\left(
-q\right) ^{a+ip^{n}}\text{.}  \label{equation 5}
\end{equation}

where%
\begin{eqnarray}
\left[ a+ip^{n}\right] _{q}^{k} &=&\sum_{j=0}^{k}\binom{k}{j}\left[ a\right]
_{q}^{k-j}q^{aj}\left[ p^{n}\right] _{q}^{j}\left[ i\right] _{q^{p^{n}}}^{j}
\label{equation 7} \\
&=&\left[ a\right] _{q}^{k}+k\left[ a\right] _{q}^{k-1}q^{a}\left[ p^{n}%
\right] _{q}\left[ i\right] _{q^{p^{n}}}+...+q^{ak}\left[ p^{n}\right]
_{q}^{k}\left[ i\right] _{q^{p^{n}}}^{k}\text{.}  \notag
\end{eqnarray}

and%
\begin{equation}
w^{a+ip^{n}}=w^{a}\sum_{l=0}^{ip^{n}}\binom{ip^{n}}{l}\left( w-1\right)
^{l}\equiv w^{a}\left( \func{mod}p^{n}\right) \text{.}  \label{equation 8}
\end{equation}

Similarly,%
\begin{equation}
\left( -q\right) ^{a+ip^{n}}=\left( -q\right) ^{a}\sum_{l=0}^{ip^{n}}\binom{%
ip^{n}}{l}\left( -1\right) ^{l}\left( q+1\right) ^{l}\equiv \left( -q\right)
^{a}\left( \func{mod}p^{n}\right) \text{.}  \label{equation 9}
\end{equation}

By (\ref{equation 5}), (\ref{equation 7}), (\ref{equation 8}) and (\ref%
{equation 9}), we have the following%
\begin{eqnarray*}
\mu _{\mathcal{P},-q}^{\left( \omega \right) }\left( a+p^{n}\mathcal{%
\mathbb{Z}
}_{p}\right) &\equiv &\left( -1\right) ^{a}\omega ^{a}q^{a}\left[ a\right]
_{q}^{k}\left( \func{mod}p^{n}\right) \\
&\equiv &\left( -1\right) ^{a}\omega ^{a}q^{a}\mathcal{P}\left( a\right)
\left( \func{mod}p^{n}\right) \text{.}
\end{eqnarray*}

For $x\in \mathcal{%
\mathbb{Z}
}_{p}$, let $x\equiv x_{n}\left( \func{mod}p^{n}\right) $ and $x\equiv
x_{n+1}\left( \func{mod}p^{n+1}\right) $, where $x_{n}$, $x_{n+1}\in 
\mathcal{%
\mathbb{Z}
}$ with $0\leq x_{n}<p^{n}$ and $0\leq x_{n+1}<p^{n+1}$.

Then, we procure the following%
\begin{equation*}
\left\vert \left[ p^{n}\right] _{-q}\mu _{\mathcal{P},-q}^{\left( \omega
\right) }\left( a+p^{n}\mathcal{%
\mathbb{Z}
}_{p}\right) -\left[ p^{n+1}\right] _{-q}\mu _{\mathcal{P},-q}^{\left(
\omega \right) }\left( a+p^{n+1}\mathcal{%
\mathbb{Z}
}_{p}\right) \right\vert \leq Cp^{-n}\text{,}
\end{equation*}

where $C$ is positive constant and $n>>0$.

Let $UD\left( \mathcal{%
\mathbb{Z}
}_{p},\mathcal{%
\mathbb{C}
}_{p}\right) $ be the space of uniformly differentiable functions on $%
\mathcal{%
\mathbb{Z}
}_{p}$ with supnorm 
\begin{equation*}
\left\Vert f\right\Vert _{\infty }=\underset{x\in 
\mathbb{Z}
_{p}}{\sup }\left\vert f\left( x\right) \right\vert _{p}.
\end{equation*}

The difference quotient $\Delta _{1}f$ of $f$ is the function of two
variables given by 
\begin{equation*}
\Delta _{1}f\left( m,x\right) =\frac{f\left( x+m\right) -f\left( x\right) }{m%
},\text{ for all }x\text{, }m\in 
\mathbb{Z}
_{p}\text{, }m\neq 0\text{.}
\end{equation*}

A function $f:%
\mathbb{Z}
_{p}\rightarrow 
\mathbb{C}
_{p}$ is said to be a Lipschitz function if there exists a constant $M>0$ $%
\left( \text{the Lipschitz constant of }f\right) $ such that%
\begin{equation*}
\left\vert \Delta _{1}f\left( m,x\right) \right\vert \leq M\text{ for all }%
m\in 
\mathbb{Z}
_{p}\backslash \left\{ 0\right\} \text{ and }x\in 
\mathbb{Z}
_{p}.
\end{equation*}

The $%
\mathbb{C}
_{p}$ linear space consisting of all Lipschitz function is denoted by $%
Lip\left( 
\mathbb{Z}
_{p},%
\mathbb{C}
_{p}\right) $. This space is a Banach space with the respect to the norm $%
\left\Vert f\right\Vert _{1}=\left\Vert f\right\Vert _{\infty }\tbigvee
\left\Vert \Delta _{1}f\right\Vert _{\infty }$ (for more informations, see 
\cite{Kim 1}, \cite{Kim 2}, \cite{Kim 3}, \cite{Kim 4}, \cite{Kim 5}, \cite%
{Kim 6}, \cite{Jang}). The objective of this paper is to introduce weighted $%
q$-Hardy Littlewood type maximal operator on the fermionic $p$-adic $q$%
-integral on $%
\mathbb{Z}
_{p}$. Also, we show that the boundedness of the weighted $q$%
-Hardy-littlewood-type maximal operator in the $p$-adic integer ring.

\section{\qquad \textbf{The weighted }$q$\textbf{-Hardy-littlewood-type
maximal operator}}

In view of (\ref{equation 2}) and the definition of fermionic $p$-adic $q$%
-integral on $%
\mathbb{Z}
_{p}$, we now consider the following theorem.

\begin{theorem}
Let $\mu _{-q}^{\left( w\right) }$ be a strongly fermionic $p$-adic $q$%
-invariant on $%
\mathbb{Z}
_{p}$ and $f\in UD\left( 
\mathbb{Z}
_{p},%
\mathbb{C}
_{p}\right) $. Then for any $n\in 
\mathbb{Z}
$ and any $\xi \in 
\mathbb{Z}
_{p}$, we have
\end{theorem}

$(1)$ $\int_{a+p^{n}%
\mathbb{Z}
_{p}}\omega ^{\xi }f\left( \xi \right) \left( -q\right) ^{-\xi }d\mu
_{-q}\left( \xi \right) =\frac{\left( -1\right) ^{a}\omega ^{a}}{\left[ p^{n}%
\right] _{-q}}\int_{%
\mathbb{Z}
_{p}}\omega ^{\xi }f\left( a+p^{n}\xi \right) \left( -q\right) ^{-p^{n}\xi
}d\mu _{-q^{p^{n}}}\left( \xi \right) $,

$(2)$ $\int_{a+p^{n}%
\mathbb{Z}
_{p}}\omega ^{\xi }d\mu _{-q}\left( \xi \right) =\frac{\omega ^{a}\left(
-q\right) ^{a}}{\left[ p^{n}\right] _{-q}}\frac{2}{1+\omega ^{p^{n}}q^{p^{n}}%
}$.

\begin{proof}
(1) By using (\ref{equation 6}) and (\ref{equation 2}), we see the
followings applications%
\begin{eqnarray*}
&&\int_{a+p^{n}%
\mathbb{Z}
_{p}}\omega ^{\xi }f\left( \xi \right) \left( -q\right) ^{-\xi }d\mu
_{-q}\left( \xi \right) \\
&=&\lim_{m\rightarrow \infty }\frac{1}{\left[ p^{m+n}\right] _{-q}}\sum_{\xi
=0}^{p^{m}-1}\omega ^{a+p^{n}\xi }f\left( a+p^{n}\xi \right) \left(
-q\right) ^{-\left( a+p^{n}\xi \right) }q^{a+p^{n}\xi }\left( -1\right)
^{a+p^{n}\xi } \\
&=&\left( -1\right) ^{a}\omega ^{a}\lim_{m\rightarrow \infty }\frac{1}{\left[
p^{m}\right] _{-q^{p^{n}}}\left[ p^{n}\right] _{-q}}\sum_{\xi
=0}^{p^{m}-1}\omega ^{\xi }\left( -q\right) ^{-p^{n}\xi }f\left( a+p^{n}\xi
\right) \left( -q^{p^{n}}\right) ^{\xi } \\
&=&\frac{\left( -1\right) ^{a}\omega ^{a}}{\left[ p^{n}\right] _{-q}}\int_{%
\mathbb{Z}
_{p}}\omega ^{\xi }f\left( a+p^{n}\xi \right) \left( -q\right) ^{-p^{n}\xi
}d\mu _{-q^{p^{n}}}\left( \xi \right) .
\end{eqnarray*}

(2) By the same method of (1), then, we easily derive the following%
\begin{eqnarray*}
&&\int_{a+p^{n}%
\mathbb{Z}
_{p}}\omega ^{\xi }d\mu _{-q}\left( \xi \right) \\
&=&\lim_{m\rightarrow \infty }\frac{1}{\left[ p^{m+n}\right] _{-q}}\sum_{\xi
=0}^{p^{m}-1}\omega ^{a+\xi p^{n}}\left( -q\right) ^{a+\xi p^{n}} \\
&=&\frac{\omega ^{a}\left( -q\right) ^{a}}{\left[ p^{n}\right] _{-q}}%
\lim_{m\rightarrow \infty }\frac{1}{\left[ p^{m}\right] _{-q^{p^{n}}}}%
\sum_{\xi =0}^{p^{m}-1}\left( \omega ^{p^{n}}\right) ^{\xi }\left(
-q^{p^{n}}\right) ^{\xi } \\
&=&\frac{\omega ^{a}\left( -q\right) ^{a}}{\left[ p^{n}\right] _{-q}}%
\lim_{m\rightarrow \infty }\frac{1+\left( \omega ^{p^{n}}q^{p^{n}}\right)
^{p^{m}}}{1+\omega ^{p^{n}}q^{p^{n}}} \\
&=&\frac{\omega ^{a}\left( -q\right) ^{a}}{\left[ p^{n}\right] _{-q}}\frac{2%
}{1+\omega ^{p^{n}}q^{p^{n}}}
\end{eqnarray*}

Since $\underset{m\rightarrow \infty }{\lim }q^{p^{m}}=1$ for $\left\vert
1-q\right\vert _{p}<1,$ our assertion follows.
\end{proof}

We are now ready to introduce definition of weighted $q$%
-Hardy-littlewood-type maximal operator related to fermionic $p$-adic $q$%
-integral on $%
\mathbb{Z}
_{p}$ with a strong fermionic $p$-adic $q$-invariant distribution $\mu _{-q}$
in the $p$-adic integer ring.

\begin{definition}
Let $\mu _{-q}^{\left( \omega \right) }$ be a strongly fermionic $p$-adic $q$%
-invariant distribution on $%
\mathbb{Z}
_{p}$ and $f\in UD\left( 
\mathbb{Z}
_{p},%
\mathbb{C}
_{p}\right) $. Then, $q$-Hardy-littlewood-type maximal operator with weight
related to fermionic $p$-adic $q$-integral on $a+p^{n}%
\mathbb{Z}
_{p}$ is defined by the following 
\begin{equation*}
\mathcal{M}_{p,q}^{\left( \omega \right) }f\left( a\right) =\underset{n\in 
\mathbb{Z}
}{\sup }\frac{1}{\mu _{1,-q}^{\left( w\right) }\left( \xi +p^{n}%
\mathbb{Z}
_{p}\right) }\int_{a+p^{n}%
\mathbb{Z}
_{p}}\omega ^{\xi }\left( -q\right) ^{-\xi }f\left( \xi \right) d\mu
_{-q}\left( \xi \right)
\end{equation*}%
for all $a\in 
\mathbb{Z}
_{p}$.
\end{definition}

We recall that famous Hardy-littlewood maximal operator $\mathcal{M}_{\mu }$%
, which is defined by 
\begin{equation}
\mathcal{M}_{\mu }f\left( a\right) =\underset{a\in Q}{\sup }\frac{1}{\mu
\left( Q\right) }\int_{Q}\left\vert f\left( x\right) \right\vert d\mu \left(
x\right) \text{,}  \label{equation 3}
\end{equation}

where $f:%
\mathbb{R}
^{k}\rightarrow 
\mathbb{R}
^{k}$ is a locally bounded Lebesgue measurable function, $\mu $ is a
Lebesgue measure on $\left( -\infty ,\infty \right) $ and the supremum is
taken over all cubes $Q$ which are parallel to the coordinate axes. Note
that the boundedness of the Hardy-Littlewood maximal operator serves as one
of the most important tools used in the investigation of the properties of
variable exponent spaces (see \cite{Jang}). The essential aim of Theorem 1
is to deal with the weighted $q$-extension of the classical Hardy-Littlewood
maximal operator in the space of $p$-adic Lipschitz functions on $%
\mathbb{Z}
_{p}$ and to find the boundedness of them. By the meaning of Definition 1,
then, we state the following theorem.

\begin{theorem}
Let $f\in UD\left( 
\mathbb{Z}
_{p},%
\mathbb{C}
_{p}\right) $ and $x\in 
\mathbb{Z}
_{p}$, we get
\end{theorem}

\textit{(1)} $\mathcal{M}_{p,q}^{\left( \omega \right) }f\left( a\right) =%
\frac{\left( -1\right) ^{a}}{2q^{a}}\underset{n\in 
\mathbb{Z}
}{\sup \left( 1+\omega ^{p^{n}q^{p^{n}}}\right) }\int_{%
\mathbb{Z}
_{p}}\omega ^{\xi }f\left( x+p^{n}\xi \right) \left( -q\right) ^{-p^{n}\xi
}d\mu _{-q^{p^{n}}}\left( \xi \right) $,

\textit{(2)} $\left\vert \mathcal{M}_{p,q}^{\left( \omega \right) }f\left(
a\right) \right\vert _{p}\leq \left\vert \frac{\left( -1\right) ^{a}}{2q^{a}}%
\right\vert _{p}\underset{n\in 
\mathbb{Z}
}{\sup }\left\vert 1+\omega ^{p^{n}}q^{p^{n}}\right\vert _{p}\left\Vert
f\right\Vert _{1}\left\Vert \left( \frac{-q^{p^{n}}}{\omega }\right)
^{-\left( .\right) }\right\Vert _{L^{1}}$,

\textit{where} $\left\Vert \left( \frac{-q^{p^{n}}}{\omega }\right)
^{-\left( .\right) }\right\Vert _{L^{1}}=\int_{%
\mathbb{Z}
_{p}}\left( \frac{-q^{p^{n}}}{\omega }\right) ^{-\xi }d\mu
_{-q^{p^{n}}}\left( \xi \right) $.

\begin{proof}
(1) Because of Theorem 1 and Definition 1, we see%
\begin{eqnarray*}
M_{p,q}^{\left( \omega \right) }f\left( a\right) &=&\underset{n\in 
\mathbb{Z}
}{\sup }\frac{1}{\mu _{1,-q}^{\left( \omega \right) }\left( \xi +p^{n}%
\mathbb{Z}
_{p}\right) }\int_{a+p^{n}%
\mathbb{Z}
_{p}}\omega ^{\xi }\left( -q\right) ^{-\xi }f\left( \xi \right) d\mu
_{-q}\left( \xi \right) \\
&=&\frac{\left( -1\right) ^{a}}{2q^{a}}\underset{n\in 
\mathbb{Z}
}{\sup \left( 1+\omega ^{p^{n}q^{p^{n}}}\right) }\int_{%
\mathbb{Z}
_{p}}\omega ^{\xi }f\left( x+p^{n}\xi \right) \left( -q\right) ^{-p^{n}\xi
}d\mu _{-q^{p^{n}}}\left( \xi \right) \text{.}
\end{eqnarray*}

(2) On account of (1), we can derive the following%
\begin{eqnarray*}
\left\vert M_{p,q}^{\left( \omega \right) }f\left( a\right) \right\vert _{p}
&=&\left\vert \frac{\left( -1\right) ^{a}}{2q^{a}}\underset{n\in 
\mathbb{Z}
}{\sup }\left( 1+\omega ^{p^{n}}q^{p^{n}}\right) \int_{%
\mathbb{Z}
_{p}}\omega ^{\xi }f\left( x+p^{n}\xi \right) \left( -q\right) ^{-p^{n}\xi
}d\mu _{-q^{p^{n}}}\left( \xi \right) \right\vert _{p} \\
&\leq &\left\vert \frac{\left( -1\right) ^{a}}{2q^{a}}\right\vert _{p}%
\underset{n\in 
\mathbb{Z}
}{\sup }\left\vert \left( 1+\omega ^{p^{n}}q^{p^{n}}\right) \int_{%
\mathbb{Z}
_{p}}\omega ^{\xi }f\left( x+p^{n}\xi \right) \left( -q\right) ^{-p^{n}\xi
}d\mu _{-q^{p^{n}}}\left( \xi \right) \right\vert _{p} \\
&\leq &\left\vert \frac{\left( -1\right) ^{a}}{2q^{a}}\right\vert _{p}%
\underset{n\in 
\mathbb{Z}
}{\sup }\left\vert 1+\omega ^{p^{n}}q^{p^{n}}\right\vert _{p}\int_{%
\mathbb{Z}
_{p}}\left\vert f\left( a+p^{n}\xi \right) \right\vert _{p}\left\vert \left( 
\frac{-q^{p^{n}}}{\omega }\right) ^{-\xi }\right\vert _{p}d\mu
_{-q^{p^{n}}}\left( \xi \right) \\
&\leq &\left\vert \frac{\left( -1\right) ^{a}}{2q^{a}}\right\vert _{p}%
\underset{n\in 
\mathbb{Z}
}{\sup }\left\vert 1+\omega ^{p^{n}}q^{p^{n}}\right\vert _{p}\left\Vert
f\right\Vert _{1}\int_{%
\mathbb{Z}
_{p}}\left\vert \left( \frac{-q^{p^{n}}}{\omega }\right) ^{-\xi }\right\vert
_{p}d\mu _{-q^{p^{n}}}\left( \xi \right) \\
&=&\left\vert \frac{\left( -1\right) ^{a}}{2q^{a}}\right\vert _{p}\underset{%
n\in 
\mathbb{Z}
}{\sup }\left\vert 1+\omega ^{p^{n}}q^{p^{n}}\right\vert _{p}\left\Vert
f\right\Vert _{1}\left\Vert \left( \frac{-q^{p^{n}}}{\omega }\right)
^{-\left( .\right) }\right\Vert _{L^{1}}\text{.}
\end{eqnarray*}

Thus, we complete the proof of theorem.
\end{proof}

We note that Theorem 2 (2) shows the supnorm-inequality for the $q$%
-Hardy-Littlewood-type maximal operator with weight on $%
\mathbb{Z}
_{p}$, on the other hand, Theorem 2 (2) shows the following inequality%
\begin{equation}
\left\Vert \mathcal{M}_{p,q}^{\left( \omega \right) }f\right\Vert _{\infty }=%
\underset{x\in 
\mathbb{Z}
_{p}}{\sup }\left\vert \mathcal{M}_{p,q}^{\left( \omega \right) }f\left(
x\right) \right\vert _{p}\leq \mathcal{K}\left\Vert f\right\Vert
_{1}\left\Vert \left( \frac{-q^{p^{n}}}{\omega }\right) ^{-\left( .\right)
}\right\Vert _{L^{1}}  \label{equation 4}
\end{equation}

where $\mathcal{K}=\left\vert \frac{\left( -1\right) ^{a}}{2q^{a}}%
\right\vert _{p}\underset{n\in 
\mathbb{Z}
}{\sup }\left\vert 1+\omega ^{p^{n}}q^{p^{n}}\right\vert _{p}$. By the
equation (\ref{equation 4}), we get the following Corollary, which is the
boundedness for weighted $q$-Hardy-Littlewood-type maximal operator with
weight on $%
\mathbb{Z}
_{p}$.

\begin{corollary}
$\mathcal{M}_{p,q}^{\left( \omega \right) }$ is a bounded operator from $%
UD\left( 
\mathbb{Z}
_{p},%
\mathbb{C}
_{p}\right) $ into $L^{\infty }\left( 
\mathbb{Z}
_{p},%
\mathbb{C}
_{p}\right) $, where $L^{\infty }\left( 
\mathbb{Z}
_{p},%
\mathbb{C}
_{p}\right) $ is the space of all $p$-adic supnorm-bounded functions with
the 
\begin{equation*}
\left\Vert f\right\Vert _{\infty }=\underset{x\in 
\mathbb{Z}
_{p}}{\sup }\left\vert f\left( x\right) \right\vert _{p}\text{,}
\end{equation*}%
for all $f\in L^{\infty }\left( 
\mathbb{Z}
_{p},%
\mathbb{C}
_{p}\right) $.
\end{corollary}

\end{document}